\documentclass[10pt]{article}
\usepackage{amsfonts}
\usepackage{amsmath}
\usepackage{mathrsfs}
\usepackage{mathrsfs,amscd,amssymb,amsthm,amsmath,bm,graphicx,psfrag,subfigure}

\makeatletter

\renewcommand{\@seccntformat}[1]{{\csname the#1\endcsname}{\normalsize .}\hspace{.5em}}
\makeatother

\usepackage{indentfirst}

\def \[{\begin{equation}}
\def \]{\end{equation}}

\def \ec{{\rm ec}}
\def \diam{{\rm diam}}
\def \ree{{\rm ree}}
\newtheorem{thm}{Theorem}[section]

\newtheorem{lem}[thm]{Lemma}
\newtheorem{cor}[thm]{Corollary}

\newenvironment{wst}
{\setlength{\leftmargini}{1.5\parindent}
 \begin{itemize}
 \setlength{\itemsep}{-1.1mm}}
{\end{itemize}}

\begin{document}
\setlength{\baselineskip}{15pt}
\begin{center}{\Large \bf On the extremal total reciprocal edge-eccentricity of trees\footnote{Financially supported by the National Natural Science Foundation of China (Grant Nos. 11271149, 11371162) and the Program for New Century Excellent Talents in University (Grant No. NCET-13-0817).}}

\vspace{4mm}

{\large Shuchao Li\footnote{E-mail: lscmath@mail.ccnu.edu.cn (S.C.
Li), zhaolfmath@163.com (L. Zhao)},\ \ Lifang Zhao}\vspace{2mm}

Faculty of Mathematics and Statistics,  Central China Normal
University, Wuhan 430079, P.R. China
\end{center}

\noindent {\bf Abstract}: The total reciprocal edge-eccentricity is a novel graph invariant with vast potential in structure activity/property relationships. This graph invariant displays high discriminating power with respect to both biological activity and physical properties. If $G=(V_G,E_G)$ is a simple connected graph, then the total reciprocal edge-eccentricity (REE) of $G$ is defined as $\xi^{ee}(G)=\sum_{uv\in E_G}(1/\varepsilon_G(u)+1/\varepsilon_G(v))$, where $\varepsilon_G(v)$ is the eccentricity of the vertex $v$. In this paper we first introduced four edge-grafting transformations to study the mathematical properties of the reciprocal edge-eccentricity of $G$. Using these elegant mathematical properties, we characterize the extremal graphs among $n$-vertex trees with given graphic parameters, such as pendants, matching number, domination number, diameter, vertex bipartition, et al. Some sharp bounds on the reciprocal edge-eccentricity of trees are determined.

\vspace{2mm} \noindent{\it Keywords}: Edge-eccentricity; Reciprocal edge-eccentricity; Matching number; Diameter; Connectivity

\vspace{2mm}


\section{\normalsize Introduction}\setcounter{equation}{0}
Throughout this paper, we only consider simple connected graph $G=(V_G, E_G)$ on $n$ vertices and $m$ edges (so $n= |V_G|$ is its order, and
$m=|E_G|$ is its size).  The \textit{distance} between two vertices $u, v$ of $G$, written $d_G(u, v),$ is the length of a shortest $u$-$v$ path in $G.$ The \textit{eccentricity} $\varepsilon_G(v)$ of a vertex $v$ is the distance between $v$ and a furthest vertex from $v$ in $G$. For any edge $e=uv\in E_G,$ we may define \textit{edge-eccentricity} of $e$ as $\ec(e)=\varepsilon_G(u)+\varepsilon_G(v);$ whereas its \textit{reciprocal edge-eccentricity} is defined as $\ree(e)=\frac{1}{\varepsilon_G(u)}+\frac{1}{\varepsilon_G(v)}.$ When the graph is clear from the context, we will omit the subscript $G$ from the notation. We follow the notation and terminology in \cite{D-I} except if otherwise stated.

Molecular descriptors play an important role in mathematical chemistry, especially in the QSPR and QSAR modeling \cite{Y-A-T1}. Among them, a special place is reserved for the so called topological indices, or graph invariant. The well-studied distance-based graph invariant probably is the \textit{Wiener index} \cite{H-W}, one of the most well used chemical indices that correlate a chemical compound's structure with the compound's physical-chemical properties. The Wiener index, introduced
in 1947, is defined as the sum of distances between all pairs of vertices, namely that
$$
  W(G) =\sum_{\{u,v\}\subseteq V_G}d_G(u, v).
$$
For more results on Wiener index one may be referred to those in  \cite{5,7,9,8}  and the references cited therein.

Another distance-based graph invariant, defined \cite{17,19} in a fully analogous manner to Wiener index, is the \textit{Harary index}, which is equal to the sum of reciprocal distances over all unordered vertex pairs in $G$, that is,
$$
   H(G)=\sum_{\{u,v\}\subseteq V_G}\frac{1}{d_G(u,v)}.
$$
For more results on Harary index, one may be referred to \cite{6,004,32,18,19,X-K-X}.

More recently, the distance-based graph invariants involving eccentricity have attracted much attention. These graph invariants mainly include the average eccentricity \cite{1}, the superaugmented eccentric connectivity index \cite{V6}, the reformed eccentric connectivity index \cite{V15}, the eccentric distance sum \cite{V18}, augmented eccentric connectivity index \cite{V19}, etc. In particular, the \textit{average eccentricity} \cite{1,003,VA00,VA}, and the \textit{eccentric distance sum} \cite{5} of the graph $G$, written by $\xi(G)$ and $\xi^d(G)$ are defined, respectively, as
$$
  \xi(G)=\frac{1}{n}\sum_{u\in V_G}\varepsilon_G(u),\ \ \ \ \ \ \ \xi^d(G)=\sum_{\{u,v\}\subseteq V_G}(\varepsilon_G(u)+\varepsilon_G(v))d_G(u,v).
$$
Recently, mathematical properties of the eccentric distance sum of graphs have been studied. Mukungunugwa and Mukwembi \cite{VM} obtained the  asymptotically sharp upper bounds on $\xi^d(G)$ with respect to the order and minimal degree of $G.$ Geng, Zhang and one of the present authors \cite{VMX} studied the relationship between $\xi^d$ and some other parameters, such as domination number, pendants and so on of trees. For more results on $\xi^d(G)$, one may be referred to \cite{005,000,VMS,MLY} and references therein.

The \textit{total edge-eccentricity} of a graph $G$ is defined as
$$
\xi^c(G)=\sum_{e=uv\in E_G}(\varepsilon_G(u)+\varepsilon_G(v)).
$$
The total edge-eccentricity of the graph $G$ can be defined alternatively as
$$
\xi^c(G)=\sum_{u\in V_G}\varepsilon_G(u) d_G(u).
$$
This graph invariant is just the \textit{eccentric connectivity index}, which is a distance-based molecular structure descriptor, proposed by Sharma, Goswami and Madan \cite{V-R-A} in 1997. The index $\xi^c(G)$ was successfully used for mathematical models of biological activities of diverse nature \cite{V6,4}. For the study of its mathematical properties one may be referred to \cite{VH,VA,VM1} and the references therein.

Bearing in mind that the relation between Wiener index and Harary index, we study here a novel graph invariant named the \textit{total reciprocal edge-eccentricity} (REE), i.e.,
$$
\xi^{ee}(G)=\sum_{e=uv\in E_G}\left(\frac{1}{\varepsilon_G(u)}+\frac{1}{\varepsilon_G(v)}\right),
$$
which can be defined alternatively as
$$
\xi^{ee}(G)=\sum_{u\in V_G}\frac{ d_G(u)}{\varepsilon_G(u)}.
$$
This graph invariant is just the \textit{the connective eccentricity index} \cite{4}. Gupta et al. \cite{4} used nonpeptide N-benzylimidazole derivatives to investigate the predictability of the total reciprocal edge-eccentricity with respect to antihypertensive activity. They showed that the results obtained using the total reciprocal edge-eccentricity were better than the corresponding values obtained using Balaban's mean square distance index and the accuracy of prediction was found to be about 80\% in the active range. Recently, Yu et al. \cite{11,VY} studied the mathematical properties of REE of trees, unicyclic graphs and general connected graphs, respectively.

Motivated from \cite{11,VY}, we mainly study the mathematical properties of total reciprocal edge-eccentricity under some edge-grafting transformations. Furthermore, extremal values of total reciprocal edge-eccentricity are also studied for some interesting classes of trees. We organize this paper as follows. In Section 3, we introduce general graph transformations that increase the total reciprocal edge-eccentricity for connected graphs.
In Section 4, sharp upper bound is established on the maximum total reciprocal edge-eccentricity of $n$-vertex trees with $k$ pendants (resp. matching number, domination number, diameter, given bipartition). The corresponding extremal graphs are identified, respectively.
\section{\normalsize Definitions and some preliminary results}
Let $G$ be a graph with $v\in V_G, uv\in E_G.$ Then $G-v$, $G-uv$ denote the graph obtained from $G$ by deleting vertex $v\in V_G$ or edge $uv\in E_G$, respectively, and this notation is naturally extended if more than one vertex or edge is deleted. $G+uv$ is obtained from $G$ by adding an edge $uv\notin E_G$.  The symbol $\sim$ denotes that two vertices in question are adjacent.

We let $N(v)$ denote the set of all the adjacent vertices of $v$ in $G$. Then let $d(v)$ denote the degree of $v$, which is defined as $d(v)=|N(v)|$. The \textit{diameter} of $G$ is the maximal distance between any two vertices of $G$. Denote by $P_n$ and $S_n$, the path and star graph on $n$ vertices, respectively. A \textit{pendant path} at $v$ in a graph $G$ is a path in which no vertex other than $v$ is incident with any edge of $G$ outside the path, where $d_G(v)$ is at least three. We call $u$ is a \textit{pendant vertex} of $G$ or a \textit{leaf} if $d(u)=1$. Let $PV(G)$ denote the set of all pendant vertices of $G.$

The vertex with the minimum eccentricity is called the \textit{center vertex} of graph $G$. Let $w$ be a center vertex of $G$ and for two vertices $x, y$ in $G$, we call $x$ a \textit{parent} of $y$ if $x\sim y$ and $d_G(w,y)=d_G(w,x)+1.$

Two distinct edges in a graph $G$ are \textit{independent} if they do not have a common end vertex. A set of pairwise independent edges of $G$ is called a \textit{matching} in $G$, while a matching of maximum cardinality is a maximum matching in $G$. The \textit{matching number} $\beta(G)$ of $G$ is the
cardinality of its maximum matching. Let $M$ be a matching of $G$, then if a vertex is incident to an edge of $M$, then it is \textit{$M$-matched}, otherwise the vertex is \textit{$M$-unmatched}. A vertex is said to be \textit{perfectly matched} if it is matched in all maximum matching of $G$.

A \textit{dominating set} $D$ of a graph $G$ is a set of vertices such that for any vertex $x$ in $G$ we have $x\in D$ or $x$  is adjacent to a vertex of $D.$ A \textit{domination number}, denoted by $\gamma$, is the minimum of the cardinalities of all dominating sets. Let $\mathscr{D}_n^\gamma$ be the class of all connected graphs of order $n$ with domination number $\gamma.$ An \textit{independent set} of a graph $G$ is a set of vertices such that any two distinct vertices of the set are not adjacent.

Given a connected graph $G$ on $n$ vertices, its vertex set can be partitioned into two subsets $V_1$ and $V_2$ such that each edge joins a vertex in $V_1$ with a vertex in $V_2$. Suppose that $V_1$ has $p$ vertices and $V_2$ has $q$ vertices with $p+q=n$. Then we say that $G$ has a $(p, q)$-\textit{bipartition} ($p\leqslant q$). 
\begin{thm}[\cite{VY}]
Let $G$ be a connected graph with diameter $d$ and the diametric path is $P=v_0v_1\ldots v_d$. Assume that $w$ is a vertex outside $P$ having pendant neighbors $w_1, w_2, \ldots, w_t$. Let
$$
   G'=G-\{ww_1, ww_2, \ldots, ww_t\}+\{v_1w_1, v_1w_2, \ldots, v_1w_t\},
$$
then we have $\xi^{ee}(G)\geqslant\xi^{ee}(G')$.
\end{thm}

\section{\normalsize Transformations }\setcounter{equation}{0}
In this section, we are to introduce four edge-drafting transformations on connected graphs. We mainly study the effect of each of these transformations on the total reciprocal edge-eccentricity.
\subsection{\normalsize $\rho$-transformation}
Let $G$ be an $n$-vertex connected graph as depicted in Fig. 1, where $wv$ is a cut edge of $G$ and $H_1,H_2$ are two non-trivial connected subgraphs with $\varepsilon_{H_1}(w)=k\geqslant l+1$. Let $G'=G-\{vx: x\in V_{H_2}(v)\}+\{wx: x\in V_{H_2}(v)\};$ see Fig. 1.  We say that $G'$ is obtained from $G$ by $\rho$-\textit{transformation}. In particular, if $H_1$ (resp. $H_2$) is a tree, Ili\'c \cite{A-I-C} used the $\rho$-transformation to study the Laplacian coefficients of trees; Geng, Zhang and one of the present authors \cite{VMX} used the $\rho$-transformation to study the eccentric distance sum of trees; Meng and one of the present authors \cite{7} used the $\rho$-transformation to study the property of additively weighted Harary index of trees. Here we show that $\rho$-transformation increases the $\xi^{ee}(G).$
\begin{figure}[h!]
\begin{center}
\psfrag{s}{$H_1$}\psfrag{r}{$w_2$}\psfrag{h}{$H_2$}
\psfrag{4}{$v$}\psfrag{3}{$w$}\psfrag{2}{$w_1$}\psfrag{1}{$w_2$}\psfrag{8}{$w_k$}\psfrag{5}{$u_1$}\psfrag{t}{$P_{l+1}$}
\psfrag{7}{$u_l$}\psfrag{b}{$G$}\psfrag{a}{$G'$}\psfrag{e}{$\Rightarrow$}
\includegraphics[width=140mm]{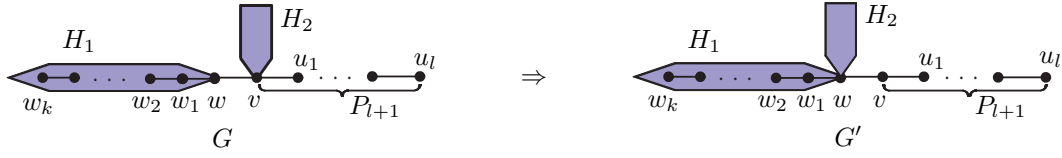}\\
\caption{$G'$ is obtained from $G$ by $\rho$-transformation.}
\end{center}
\end{figure}
\begin{thm}
Let $G'$ be the graph obtained from $G$ by $\rho$-transformation; see Fig. 1. Then $\xi^{ee}(G)<\xi^{ee}(G')$.
\end{thm}
\begin{proof}
Assume that $\varepsilon_{H_2}(v)=q$, $\varepsilon_{H_1}(w)=d_{H_1}(w, w_k)$ and $P_{k+1}=ww_1w_2\ldots w_k$ is the path between $w$ and a furthest vertex $w_k$ from $w$ in $H_1$. 

For every vertex $x\in V_{H_1}\setminus V_{P_{k+1}}$, it is routine to check that $d_G(x)=d_{G'}(x)$, $\varepsilon_G(x)=\max\{\varepsilon_{H_1}(x), d_{H_1}(x, w)+1+l, d_{H_1}(x, w)+1+q\}$, $\varepsilon_{G'}(x)=\max\{\varepsilon_{H_1}(x), d_{H_1}(x, w)+1+l, d_{H_1}(x, w)+q\}$, i.e., $\varepsilon_G(x)\geqslant \varepsilon_{G'}(x)$. Hence,
$$
   \sum_{x\in V_{H_1}\setminus V_{P_{k+1}}}\left(\frac{d_G(x)}{\varepsilon_G(x)}-\frac{d_{G'}(x)}{\varepsilon_{G'}(x)}\right)\leqslant0.
$$

For every vertex $x\in V_{H_2}\setminus \{v\}$, we have $d_G(x)=d_{G'}(x)$, $\varepsilon_G(x)=\max\{\varepsilon_{H_2}(x), d_{H_2}(x, v)+1+k\}$, $\varepsilon_{G'}(x)=\max\{\varepsilon_{H_2}(x), d_{H_2}(x, v)+k\}$, i.e., $\varepsilon_G(x)\geqslant \varepsilon_{G'}(x).$ Hence,
$$
\sum_{x\in V_{H_2}\setminus \{v\}}\left(\frac{d_G(x)}{\varepsilon_G(x)}-\frac{d_{G'}(x)}{\varepsilon_{G'}(x)}\right)\leqslant0.
$$

Consider the vertex $w_i$, $i=1, 2, \ldots, k$. Since $\varepsilon_{H_1}(w)=k$, $\varepsilon_{H_1}(w_i)\leqslant i+k$, $i=1, 2, \ldots, k$. Otherwise, assume  $m$\, ($m\in [1, k]$) is the maximum subscript such that $\varepsilon_{H_1}(w_m)\geqslant m+k+1$. Suppose that $P$ is the shortest path between $w_m$ and $t$ such that $\varepsilon_{H_1}(w_m)=|E_P|$ and $p$ (possibly equal to 0) is the minimum subscript of vertex of $V_P\cap V_{P_{k+1}}$, and we denote $w=w_0$. Clearly $t\notin V_{P_{k+1}}$. If $p=m$, we claim that there exists a shortest path between $w$ and any vertex of $P$ which contains vertex $w_m$. Otherwise, combined with $\varepsilon_{H_1}(w)=k$, i.e., $d_{H_1}(w, t)\leqslant k$, we can easily get $\varepsilon_{H_1}(w_m)\leqslant m+k$, a contradiction. Thus, $d_{H_1}(w, t)=d_{H_1}(w, w_m)+d_{H_1}(w_m, t)\geqslant 2m+k+1>k$, i.e., $\varepsilon_{H_1}(w)>k$, a contradiction. Similarly, if $p<m$, then $d_{H_1}(w, t)=d_{H_1}(w, w_p)+|E_P|-d_{H_1}(w_m, w_p)\geqslant k+1+2d_{H_1}(w, w_p)>k$, a contradiction. Hence, $\varepsilon_{H_1}(w_i)\leqslant i+k$, $i=1, 2, \ldots, k$. Therefore,
\begin{equation}\label{eq:3.2}
    \varepsilon_G(w_i)=\max\{i+k, i+1+q, i+1+l\}, \ \ \ \varepsilon_{G'}(w_i)=\max\{i+k, i+q, i+1+l\}.
\end{equation}
Note also that $d_G(u_l)=d_{G'}(u_l)=1$ and for $i=1, 2, \ldots, l-1$,
\begin{eqnarray}
  d_G(u_i)=d_{G'}(u_i)=2,\ \ \varepsilon_G(u_i)=\max\{i+q, i+1+k\}, \ \ \varepsilon_{G'}(u_i)=\max\{i+1+q, i+1+k\}
\end{eqnarray}
and
\begin{align}
  & d_G(w)=d_{H_1}(w)+1, \ \ d_{G'}(w)=d_{H_1}(w)+1+d_{H_2}(v),\ \ \varepsilon_G(w)=\max\{k, 1+q\},\ \  \varepsilon_{G'}(w)=\max\{k, q\}, \\
  & d_G(v)=d_{H_2}(v)+2,\ \  d_{G'}(v)=2,\ \ \varepsilon_G(v)=\max\{1+k, q\},\ \ \varepsilon_{G'}(v)=\max\{1+k, 1+q\}.\label{eq:3.5}
\end{align}

From the definition of $\xi^{ee}$, one has
\begin{eqnarray*}
\xi^{ee}(G)-\xi^{ee}(G')&=&\sum_{x\in V_{H_1}\setminus V_{P_{k+1}}}\left(\frac{d_G(x)}{\varepsilon_G(x)}-\frac{d_{G'}(x)}{\varepsilon_{G'}(x)}\right)+\sum_{x\in V_{H_2}\setminus \{v\}}\left(\frac{d_G(x)}{\varepsilon_G(x)}-\frac{d_{G'}(x)}{\varepsilon_{G'}(x)}\right)
\end{eqnarray*}
\begin{eqnarray}
&&+\sum_{i=1}^{k}\left(\frac{d_G(w_i)}{\varepsilon_G(w_i)}-\frac{d_{G'}(w_i)}{\varepsilon_{G'}(w_i)}\right)
+\sum_{i=1}^{l}\left(\frac{d_G(u_i)}{\varepsilon_G(u_i)}-\frac{d_{G'}(u_i)}{\varepsilon_{G'}(u_i)}\right)\notag\\
&&+\frac{d_G(w)}{\varepsilon_G(w)}-\frac{d_{G'}(w)}{\varepsilon_{G'}(w)}+\frac{d_G(v)}{\varepsilon_G(v)}-\frac{d_{G'}(v)}{\varepsilon_{G'}(v)}\notag\\
&\leqslant&\sum_{i=1}^{k}\left(\frac{d_G(w_i)}{\varepsilon_G(w_i)}-\frac{d_{G'}(w_i)}{\varepsilon_{G'}(w_i)}\right)
+\sum_{i=1}^{l}\left(\frac{d_G(u_i)}{\varepsilon_G(u_i)}-\frac{d_{G'}(u_i)}{\varepsilon_{G'}(u_i)}\right)\notag\\
&&+\frac{d_G(w)}{\varepsilon_G(w)}-\frac{d_{G'}(w)}{\varepsilon_{G'}(w)}+\frac{d_G(v)}{\varepsilon_G(v)}-\frac{d_{G'}(v)}{\varepsilon_{G'}(v)}.\label{eq:3.6}
\end{eqnarray}
We proceed by considering the following possible cases.

\textbf{Case 1}.\ $q\geqslant k+1\geqslant l+2$.\ In this case, in view of (\ref{eq:3.2})-(\ref{eq:3.5}), we have $\varepsilon_G(w_i)=i+1+q$, $\varepsilon_{G'}(w_i)=i+q$ for $1\leqslant i\leqslant k$; $\varepsilon_G(u_i)=i+q$, $\varepsilon_{G'}(u_i)=i+1+q$ for $1\leqslant i\leqslant l$; $\varepsilon_G(w)=1+q$, $\varepsilon_{G'}(w)=q$; $\varepsilon_G(v)=q$, $\varepsilon_{G'}(v)=1+q$. Thus, together with (\ref{eq:3.6}), we obtain
\begin{eqnarray*}
\xi^{ee}(G)-\xi^{ee}(G')&\leqslant & \sum_{i=1}^{k}\left(\frac{d_G(w_i)}{i+1+q}-\frac{d_{G'}(w_i)}{i+q}\right)+\sum_{i=1}^{l-1}\left(\frac{2}{i+q}-\frac{2}{i+1+q}\right)+\frac{1}{l+q}-\frac{1}{l+1+q}\\
&&+\frac{d_{H_1}(w)+1}{1+q}-\frac{d_{H_1}(w)+1+d_{H_2}(v)}{q}+\frac{d_{H_2}(v)+2}{q}-\frac{2}{1+q}\\
&=&\sum_{i=1}^{l-1}\left(\frac{d_{G}(w_i)-2}{i+1+q}-\frac{d_G(w_i)-2}{i+q}\right)+\frac{d_G(w_l)-1}{l+1+q}-\frac{d_G(w_l)-1}{l+q}\\
&&+\sum_{i=l+1}^k\left(\frac{d_{G}(w_i)}{i+1+q}-\frac{d_G(w_i)}{i+q}\right)+\frac{d_{H_1}(w)-1}{1+q}-\frac{d_{H_1}(w)-1}{q}\\
&<&0.
\end{eqnarray*}
The last inequality follows from the fact that $d_G(w_i)\geqslant 2$ for $1\leqslant i\leqslant l-1$, $d_G(w_l)\geqslant 2$ for $k\geqslant l+1$ and $d_{H_1}(w)\geqslant 2$.

\textbf{Case 2}.\ $k\geqslant q\geqslant l+1$.\ In this case, in view of (\ref{eq:3.2})-(\ref{eq:3.5}), we have $\varepsilon_G(w_i)\geqslant \varepsilon_{G'}(w_i)$ for $1\leqslant i\leqslant k$, $\varepsilon_G(u_i)=\varepsilon_{G'}(u_i)=i+1+k$ for $1\leqslant i\leqslant l$, $\varepsilon_G(w)=\max\{k, 1+q\}\geqslant\varepsilon_{G'}(w)=k$ and $\varepsilon_G(v)=\varepsilon_{G'}(v)=1+k$.
Hence, together with (\ref{eq:3.6}), we obtain
\begin{eqnarray*}
\xi^{ee}(G)-\xi^{ee}(G')&\leqslant& \frac{d_{H_1}(w)+1}{\varepsilon_{G}(w)}-\frac{d_{H_1}(w)+1+d_{H_2}(v)}{\varepsilon_{G'}(w)}+\frac{d_{H_2}(v)+2}{1+k}-\frac{2}{1+k}\\
&\leqslant&\frac{-d_{H_2}(v)}{k}+\frac{d_{H_2}(v)}{1+k}\\
&<&0.
\end{eqnarray*}

\textbf{Case 3}.\ $k\geqslant l+1\geqslant q+1$. In this case, in view of (\ref{eq:3.2})-(\ref{eq:3.5}), we have $\varepsilon_G(w_i)=\varepsilon_{G'}(w_i)$ for $1\leqslant i\leqslant k$, $\varepsilon_G(u_i)=\varepsilon_{G'}(u_i)=i+1+k$ for $1\leqslant i\leqslant l$, $\varepsilon_G(w)=\varepsilon_{G'}(w)=k$ and $\varepsilon_G(v)=\varepsilon_{G'}(v)=1+k$. Combining with (\ref{eq:3.6}) yields
\begin{eqnarray*}
\xi^{ee}(G)-\xi^{ee}(G')&\leqslant&\frac{d_{H_1}(w)+1}{\varepsilon_{G}(w)}-\frac{d_{H_1}(w)+1+d_{H_2}(v)}{\varepsilon_{G'}(w)}+\frac{d_{H_2}(v)+2}{1+k}-\frac{2}{1+k}\\
&=&\frac{-d_{H_2}(v)}{k}+\frac{d_{H_2}(v)}{1+k}\\
&<&0.
\end{eqnarray*}

This completes the proof.
\end{proof}

The following corollary is a direct consequence of Theorem 3.1.
\begin{cor}
Assume that $P_{m+l}=v_0v_1\ldots v_{m-1}v_mv_{m+1}\ldots v_{m+l}$ is a path and $u$ is a vertex of a connected graph $H$. The graph $G_{l, m}$ is obtained from $P_{m+l}$ and $H$ by identifying $u$ with $v_m$, while $G_{l-1, m+1}$ is the graph obtained from $G_{l, m}$ by moving $H$ from $v_m$ to $v_{m+1}$. If $l\geqslant m+2$, then $\xi^{ee}(G_{l, m})<\xi^{ee}(G_{l-1, m+1})$.
\end{cor}

\subsection{\normalsize $\alpha$-transformation}
Let $G_1$ be a simple graph as depicted in Fig. 2, where $H_1, H_2$ are two non-trivial connected graphs. Let $G_2=G_1-\{v_lx: x\in N_{H_2}(v_l)\}+\{v_1x: x\in N_{H_2}(v_l)\}$. We call that $G_2$ is obtained by $\alpha$-\textit{transformation} on $G_1$. In particular, if $G_1$ is a tree, Kelmans \cite{A-K-K} used this tree-transformation as depicted in Fig. 2 to prove some results on the number of spanning trees of graphs in 1976. Recently, Bollob\'as and Tyomkyn \cite{B-B-M} used this tree-transformation to count the total number of walks (resp. closed walks, paths) of trees. Here we are to show that $\alpha$-transformation increases the total reciprocal edge-eccentricity of a connected graph.
\begin{figure}[h!]
\begin{center}
\psfrag{a}{$H_1$}\psfrag{f}{$\Rightarrow$}
\psfrag{b}{$H_2$}
\psfrag{c}{$v_l$}
\psfrag{d}{$P_l$}
\psfrag{1}{$v_1$}
\psfrag{2}{$v_2$}
\psfrag{g}{$H_1$}
\psfrag{A}{$G_1$}
\psfrag{B}{$G_2$}
\includegraphics[width=110mm]{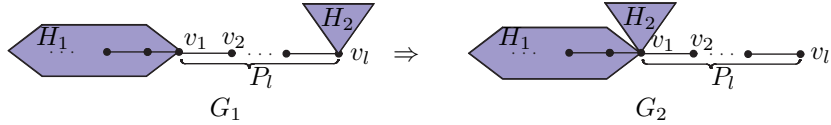}
  \caption{ $\alpha$-transformation. }
\end{center}
\end{figure}
\begin{thm}
Let $G_2$ be the connected graph that obtained from $G_1$ by $\alpha$-transformation as depicted in Fig. 2, then $\xi^{ee}(G_1)<\xi^{ee}(G_2)$.
\end{thm}
\begin{proof}
For the subgraphs $H_1$ and $H_2$, we assume, without loss of generality, that $\varepsilon_{H_1}(v_1)=q$, $\varepsilon_{H_2}(v_l)=p$ satisfying $q\geqslant p$.

For every vertex $x\in V_{H_1}\setminus \{v_1\}$, we have $d_{G_1}(x)=d_{G_2}(x)$, $\varepsilon_{G_1}(x)=\max\{\varepsilon_{H_1}(x), d_{H_1}(x, v_1)+l-1+p\}$, $\varepsilon_{G_2}(x)=\max\{\varepsilon_{H_1}(x), d_{H_1}(x, v_1)+l-1, d_{H_1}(x, v_1)+p\}$. Hence, $\varepsilon_{G_1}(x)\geqslant \varepsilon_{G_2}(x)$ for $x\in V_{H_1}\setminus \{v_1\}$.

For any vertex $x\in V_{H_2}\setminus \{v_l\}$, we have $d_{G_1}(x)=d_{G_2}(x)$, $\varepsilon_{G_1}(x)=\max\{\varepsilon_{H_2}(x),\, d_{H_2}(x, v_l)+l-1+q\}$, $\varepsilon_{G_2}(x)=\max\{\varepsilon_{H_2}(x),\, d_{H_2}(x, v_l)+q, d_{H_2}(x, v_l)+l-1\}$. Thus, $\varepsilon_{G_1}(x)\geqslant\varepsilon_{G_2}(x)$ for $x\in V_{H_2}\setminus \{v_l\}$.

For every vertex $v_i\in V_{P_l}$, $i=1, 2, \ldots, l$, we have $\varepsilon_{G_1}(v_i)=\max\{l-i+p, i-1+q\}$ and $\varepsilon_{G_2}(v_i)=\max\{l-i, i-1+q\}$. Thus we have $\varepsilon_{G_1}(v_i)\geqslant \varepsilon_{G_2}(v_i)$ for $v_i\in V_{P_l}$. What's more, $\varepsilon_{G_2}(v_l)>\varepsilon_{G_2}(v_1)$,  $d_{G_1}(v_i)=d_{G_2}(v_i)$ except for $d_{G_1}(v_1)=d_{H_1}(v_1)+1$, $d_{G_2}(v_1)=d_{H_1}(v_1)+1+d_{H_2}(v_l)$ and $d_{G_1}(v_l)=d_{H_2}(v_l)+1$, $d_{G_2}(v_l)=1$.

Therefore,
\begin{eqnarray*}
\xi^{ee}(G_1)-\xi^{ee}(G_2)&=&\sum_{x\in V_{H_1}\setminus \{v_1 \}}\left(\frac{d_{G_1}(x)}{\varepsilon_{G_1}(x)}-\frac{d_{G_2}(x)}{\varepsilon_{G_2}(x)}\right)+\sum_{x\in V_{H_2}\setminus \{ v_l\}}\left(\frac{d_{G_1}(x)}{\varepsilon_{G_1}(x)}-\frac{d_{G_2}(x)}{\varepsilon_{G_2}(x)}\right)\\
\ \ \ \ \ &&+\sum_{i=1}^{l}\left(\frac{d_{G_1}(v_i)}{\varepsilon_{G_1}(v_i)}-\frac{d_{G_2}(v_i)}{\varepsilon_{G_2}(v_i)}\right)\\
&\leqslant&\frac{d_{G_1}(v_1)}{\varepsilon_{G_1}(v_1)}-\frac{d_{G_2}(v_1)}{\varepsilon_{G_2}(v_1)}+\frac{d_{G_1}(v_l)}{\varepsilon_{G_1}(v_l)}-\frac{d_{G_2}(v_l)}{\varepsilon_{G_2}(v_l)}\\
&\leqslant&-\frac{d_{H_2}(v_l)}{\varepsilon_{G_2}(v_1)}+\frac{d_{H_2}(v_l)}{\varepsilon_{G_2}(v_l)}\\
&=&d_{H_2}(v_l)\left(\frac{1}{\varepsilon_{G_2}(v_l)}-\frac{1}{\varepsilon_{G_2}(v_1)}\right)\\
&<&0.
\end{eqnarray*}
Thus, $\xi^{ee}(G_1)<\xi^{ee}(G_2)$, as desired.
\end{proof}

The following result is a direct consequence of Theorem 3.3.
\begin{cor}[\cite{11}]
Let $H_1$ and $H_2$ be two disjoint connected graphs each of which contains at least $2$ vertices with $u\in V_{H_1}$, $v\in V_{H_2}$. Let $G_1$ be the graph obtained from $H_1\cup H_2$ by adding an edge $uv$. Let $G_2$ be the graph obtained from $H_1\cup H_2$ by identifying $u$ and $v$ (to a new vertex, say $u$) and adding a pendant edge, say $uv$ without confusion. Then $\xi^{ee}(G_1)<\xi^{ee}(G_2)$.
\end{cor}

\subsection{\normalsize $\theta$-transformation}
\begin{figure}[h!]
\begin{center}
\psfrag{a}{$w$}
\psfrag{b}{$u$}
\psfrag{c}{$v$}
\psfrag{d}{$v_1$}
\psfrag{e}{$v_2$}
\psfrag{f}{$v_k$}
\psfrag{j}{$G$}
\psfrag{1}{$H$}\psfrag{g}{$\Rightarrow$}
\psfrag{6}{$H_2$}
\psfrag{7}{$H_1$}
\psfrag{k}{$G[v\rightarrow w;2]$}
\includegraphics[width=80mm]{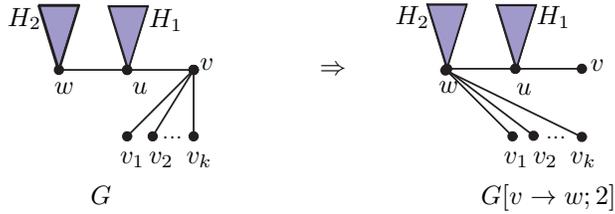}
  \caption{ $\theta$-transformation }
\end{center}
\end{figure}

Let $G$ be the connected graph as depicted in Fig. 3, where $H_2$ is a non-trivial subgraph and $uw$ is a cut edge of $G$ with $d_G(w)\geqslant2$. Let $G[v\rightarrow w;2]$ be the graph obtained from $G$ by deleting all pendant edges $vz, z\in W$ and adding all pendant edges $wz, z\in W$, where $W=N_{G}(v)\backslash\{u\}.$ In notation,
$$
   G[v\rightarrow w;2]=G-\{vz: z\in W\}+\{wz: z\in W\}
$$
and we say $G[v\rightarrow w;2]$ is obtained from $G$ by the $\theta$-\textit{transformation}. If $H_1, H_2$ are two bipartite graphs, Geng, Zhang and one of the present authors \cite{VMS} used the $\theta$-transformation to study the eccentric distance sum of trees. Li, Zhang and one of the present authors \cite{L-L-Z} used the $\theta$-transformation to study the Laplacian permanent of trees with given bipartition. Next, we are to use the $\theta$-transformation as a tool to study the total reciprocal edge-eccentricity of trees.
\begin{thm}
Let $G[v\rightarrow w;2]$ be the graph obtained from $G$ by the $\theta$-transformation (see Fig. 3). If $\varepsilon_{H_2}(w) \geqslant \varepsilon_{H_1}(u)$, then we have $\xi^{ee}(G)<\xi^{ee}(G[v\rightarrow w;2]).$
\end{thm}
\begin{proof}
For convenience, denote $G'=G[v\rightarrow w;2]$, $\varepsilon_{H_1}(u)=d_1$ and $\varepsilon_{H_2}(w)=d_2$. Since $uw$ is a cut edge of $G$ and $d_G(w)\geqslant2$, it is obviously $d_2\geqslant1$, i.e., $V_{H_2}\setminus \{w\}\neq \emptyset$.

For every vertex $x\in V_{H_1}, d_G(x)=d_{G'}(x), \varepsilon_G(x)=\max\{\varepsilon_{H_1}(x), d_G(x, u)+1+d_2\}=\varepsilon_{G'}(x)$.

For every vertex $x\in V_{H_2}\setminus \{w\}, d_G(x)=d_{G'}(x)$, $\varepsilon_G(x)=\max\{\varepsilon_{H_2}(x), d_G(x, w)+1+d_1,d_G(x, w)+3\}$, $\varepsilon_{G'}(x)=\max\{\varepsilon_{H_2}(x), d_G(x, w)+1+d_1,d_G(x, w)+2\}$. Hence, $\varepsilon_G(x)\geqslant\varepsilon_{G'}(x)$ for $x\in V_{H_2}\setminus \{w\}.$

For every vertex $v_i$, $i=1, 2, \ldots, k$, $d_G(v_i)=d_{G'}(v_i)=1$, $\varepsilon_G(v_i)=3+d_2,$ $\varepsilon_{G'}(v_i)=\max\{1+d_2, d_1+2, 3\}$, i.e., $\varepsilon_G(v_i)>\varepsilon_{G'}(v_i)$ for $d_2\geqslant d_1$, $d_2\geqslant1,\, i=1,2,\ldots, k$.

For the vertex $w$, one has $d_G(w)=d_{H_2}(w)+1$, $d_G(v)=k+1$ and $d_{G'}(w)=d_{H_2}(w)+1+k$, $d_{G'}(v)=1$. $\varepsilon_G(w)=\max\{d_2, 1+d_1, 3\}, \varepsilon_{G'}(w)=\max\{d_2, 1+d_1, 2\}$, i.e., $\varepsilon_G(w)\geqslant\varepsilon_{G'}(w).$ In view of Fig.~3, we have $\varepsilon_G(v)=\varepsilon_{G'}(v)=d_2+2$.

Therefore,
\begin{eqnarray*}
\xi^{ee}(G)-\xi^{ee}(G')&=&\sum_{x\in V_{H_1}}\left(\frac{d_G(x)}{\varepsilon_G(x)}-\frac{d_{G'}(x)}{\varepsilon_{G'}(x)}\right)+\sum_{x\in V_{H_2}\setminus\{w\}}\left(\frac{d_G(x)}{\varepsilon_G(x)}-\frac{d_{G'}(x)}{\varepsilon_{G'}(x)}\right)\\
&&+\sum_{i=1}^k\left(\frac{d_G(v_i)}{\varepsilon_G(v_i)}-\frac{d_{G'}(v_i)}{\varepsilon_{G'}(v_i)}\right)+\frac{d_G(w)}{\varepsilon_G(w)}-\frac{d_{G'}(w)}{\varepsilon_{G'}(w)}+\frac{d_G(v)}{\varepsilon_G(v)}-\frac{d_{G'}(v)}{\varepsilon_{G'}(v)}\\
&<&-\frac{k}{\varepsilon_{G'}(w)}+\frac{k}{\varepsilon_{G'}(v)}
\end{eqnarray*}
\begin{eqnarray*}
&=&k\left(\frac{1}{\varepsilon_{G'}(v)}-\frac{1}{\varepsilon_{G'}(w)}\right)\\
&<&0,
\end{eqnarray*}
The last inequality follows by $\varepsilon_{G'}(v)=d_2+2>\varepsilon_{G'}(w)=\max\{d_2, 1+d_1, 2\}$.

This completes the proof.
\end{proof}
\begin{thm}
Let $G$ be a connected graph with diameter $d$. Choose a diametral path $P=v_0v_1\ldots v_d$ of $G$ such that there exist a pendant $w_1$ not being attached on $P.$ Denote the unique neighbor of $w_1$ by $w.$ Let
$$  G_1=\left\{
             \begin{array}{ll}
              G-ww_1+v_{\frac{d}{2}-1}w_1, & \hbox{if $d$ is even,}\\[5pt]
               G-ww_1+v_{\frac{d-1}{2}-1}w_1, & \hbox{if $d$ is odd.}
             \end{array}
           \right.
$$
Then $\xi^{ee}(G)\leqslant \xi^{ee}(G_1)$, and the equality holds if and only if $\varepsilon_G(w)=\varepsilon_G(v_{\frac{d}{2}-1})=\frac{d}{2}+1$ if $d$ is even and $\varepsilon_G(w)=\varepsilon_G(v_{\frac{d-1}{2}-1})=\frac{d+1}{2}+1$ otherwise.
\end{thm}
\begin{proof}
If $d$ is even, we have $\diam(G)=\diam(G_1)$, $d_G(x)=d_{G_1}(x)$ and $\varepsilon_G(x)=\varepsilon_{G_1}(x)$ for $x\in V_G\setminus \{w,w_1, v_{\frac{d}{2}-1}\}$. $d_G(w)=d_{G_1}(w)+1$, $\varepsilon_G(w)=\varepsilon_{G_1}(w)$; $d_G(w_1)=d_{G_1}(w_1)=1,$ $\varepsilon_G(w_1)=\varepsilon_G(w)+1$, $\varepsilon_{G_1}(w_1)=\varepsilon_{G_1}(v_{\frac{d}{2}-1})+1=\varepsilon_G(v_{\frac{d}{2}-1})+1$; $d_G(v_{\frac{d}{2}-1})=d_{G_1}(v_{\frac{d}{2}-1})-1$, $\varepsilon_G(v_{\frac{d}{2}-1})=\varepsilon_{G_1}(v_{\frac{d}{2}-1})$. Hence, by the definition of the REE we have
\begin{eqnarray*}
\xi^{ee}(G)-\xi^{ee}(G_1)&=&\sum_{x\in V_G} \frac{d_G(x)}{\varepsilon_G(x)}-\sum_{x\in V_G}\frac{d_{G_1}(x)}{\varepsilon_{G_1}(x)}\\
&=&\frac{d_G(v_{\frac{d}{2}-1})}{\varepsilon_G(v_{\frac{d}{2}-1)}}-\frac{d_{G_1}(v_{\frac{d}{2}-1})}{\varepsilon_{G_1}(v_{\frac{d}{2}-1)}}+
\frac{d_G(w)}{\varepsilon_G(w)}-\frac{d_{G_1}(w)}{\varepsilon_{G_1}(w)}+\frac{d_G(w_1)}{\varepsilon_G(w_1)}-\frac{d_{G_1}(w_1)}{\varepsilon_{G_1}(w_1)}\\
&=&-\frac{1}{\varepsilon_G(v_{\frac{d}{2}-1})}+\frac{1}{\varepsilon_G(w)}+\frac{1}{\varepsilon_G(w)+1}-\frac{1}{\varepsilon_G(v_{\frac{d}{2}-1})+1}\\
&\leqslant&0,
\end{eqnarray*}
The last inequality follows by $\varepsilon_G(v_{\frac{d}{2}-1})=\frac{d}{2}+1\leqslant\varepsilon_G(w)$ with equality if and only if $\varepsilon_G(w)=\varepsilon_G(v_{\frac{d}{2}-1})=\frac{d}{2}+1$.

By a similar discussion as in the proof for the even $d$, we may also obtain, for odd $d$, $\xi^{ee}(G)\leqslant \xi^{ee}(G_1)$, and the equality holds if and only if $\varepsilon_G(w)=\varepsilon_G(v_{\frac{d-1}{2}-1})=\frac{d+1}{2}+1$. We omit the procedure here.
\end{proof}

\section{\normalsize Some applications of the edge-grafting theorems on the REE of trees }\setcounter{equation}{0}
In this section, we mainly use the four edge-grafting transformation theorems established in Section 3 to study the total reciprocal edge-eccentricity of trees with some given parameters, such as pendants, matching number, domination number, given bipartition, and so on. Some sharp bounds on REE are obtained. The corresponding extremal graphs are identified respectively.

Let $\mathscr{T}^k_n$ be the set of all $n$-vertex trees with $k$ leaves. Clearly, $\mathscr{T}^{n-1}_n=\{S_n\}$ and $\mathscr{T}^2_n=\{P_n\}.$ So in what follows, we consider $\mathscr{T}^k_n$ for $3\leqslant k\leqslant n-2.$ A \textit{spider} is a tree with at most one vertex of degree more than 2, called the \textit{hub} of the spider, otherwise any vertex can be hub. A \textit{leg} of a spider is a path from the hub to one of its leaves. Let $S(a_1, a_2, \ldots, a_k)$ be a spider with $k$ legs $L_1, L_2, \ldots, L_k$ such that $|E_{L_i}|=a_i$ ($i=1, 2, \ldots, k$) satisfying $\sum_{i=1}^ka_i=n-1$. If $|a_i-a_j|\leqslant 1$ for $1\leqslant i, j\leqslant k$, then we call $S(a_1, a_2, \ldots, a_k)$ a \textit{balanced spider}.

 Let $S'(a_1, \ldots, a_s, a_1', \ldots, a_t')$ be a tree obtained from two spiders $S(a_1, \ldots, a_s)$, $S(a_1', \ldots, a_t')$ by adding an edge $uv$, where $u$, $v$ are the center vertices of $S(a_1, \ldots, a_s)$, $S(a_1', \ldots, a_t')$, respectively, and the eccentricity of $u$ and $v$ are equal in $S'(a_1, \ldots, a_s, a_1', \ldots, a_t')$, $s, t>1$, $s+t=k$.

 Let $\mathcal{S}_n^k=\{S'(a_1, \ldots, a_s, a_1', \ldots, a_t')|a_1=\cdots=a_s=a_1'=\cdots=a_t', s, t>1, s+t=k\}$. Clearly, $n\equiv 2 \pmod{k}$ for any tree $T\in \mathcal{S}_n^k$ and $\xi^{ee}(T_1)=\xi^{ee}(T_2)$ for any $T_1,T_2\in \mathcal{S}_n^k$.

\begin{thm}
Let $T$ be in $\mathscr{T}^k_n$ with the maximal REE. Then $T$ is the balanced spider $S(\underbrace{\lfloor\frac{n-1}{k}\rfloor, \ldots, \lfloor\frac{n-1}{k}\rfloor}_{k-r},\linebreak \underbrace{\lceil\frac{n-1}{k}\rceil, \ldots, \lceil\frac{n-1}{k}\rceil}_r)$ or $T\in \mathcal{S}_n^k$ if $n\equiv 2 \pmod{k}$, and $T$ is the balanced spider $S(\underbrace{\lfloor\frac{n-1}{k}\rfloor, \ldots, \lfloor\frac{n-1}{k}\rfloor}_{k-r}, \linebreak \underbrace{\lceil\frac{n-1}{k}\rceil, \ldots, \lceil\frac{n-1}{k}\rceil}_r)$, otherwise.
\end{thm}
\begin{proof}
Choose $T$ in $\mathscr{T}^k_n$ such that $\xi^{ee}(T)$ is as large as possible. By Theorem 3.1, we get $T$ is a spider or $T\cong S'(a_1, \ldots, a_s, a_1', \ldots, a_t')$.

If $T$ is a spider, then $T$ must be a balanced spider. Otherwise, by Corollary 3.2, there exists another spider $T'$ such that $\xi^{ee}(T)<\xi^{ee}(T')$, a contradiction. That is, $T\cong S(\underbrace{\lfloor\frac{n-1}{k}\rfloor, \ldots, \lfloor\frac{n-1}{k}\rfloor}_{k-r}, \underbrace{\lceil\frac{n-1}{k}\rceil, \ldots, \lceil\frac{n-1}{k}\rceil}_r)$.

Similarly, if $T\cong S'(a_1, \ldots, a_s, a_1', \ldots, a_t')$, then by Corollary 3.2, $|a_i-a_j|\leqslant 1$ and $|a_h'-a_k'|\leqslant 1$, $1\leqslant i, j \leqslant s$, $1\leqslant h, k\leqslant t$. If $a_1=\ldots=a_s=a_1'\ldots=a_t'$, i.e., $n\equiv 2 \pmod{k}$, then $T\in \mathcal{S}_n^k$. Otherwise, assume $1\leqslant a_1\leqslant \ldots\leqslant a_s$ and $1\leqslant a_1'\leqslant \ldots\leqslant a_t'$, then we have $a_s=a_t'$ as the eccentricity of $u$ and $v$ are equal in $T$, where $u$, $v$ are the vertices defined as above. Let $\{v_1, v_2, \ldots, v_t\}=N_T(v)\setminus \{u\}$,
$$T''=T-\{vv_1,vv_2, \ldots, vv_{t-1}\}+\{uv_1, uv_2, \ldots, uv_{t-1}\}.$$

Then, for every vertex $x\in V_T$, $d_T(x)=d_{T''}(x)$ except $d_T(u)=d_{T''}(u)-t+1$, $d_T(v)=d_{T''}(v)+t-1$ and $\varepsilon_{T''}(x)=\varepsilon_T(x)$ for any vertex $x\in V_T$. Thus, by a simple calculation, we have $\xi^{ee}(T)=\xi^{ee}(T'')$. And by Corollary 3.2, there exists a balanced spider $S'$ such that $\xi^{ee}(T'')<\xi^{ee}(S')$, i.e., $\xi^{ee}(T)<\xi^{ee}(S')$, a contradiction.

In conclusion, $T\in \mathscr{T}^k_n$ has maximal REE, then $T\in \mathcal{S}^k_n$ or $T$ is a balanced spider if $n\equiv 2 \pmod{k}$, otherwise, $T$ is a balanced spider.

This completes the proof.
\end{proof}
\begin{figure}[h!]
\begin{center}
\psfrag{1}{$n-2\beta+1$}
\psfrag{2}{$\beta-1$}
\includegraphics[width=70mm]{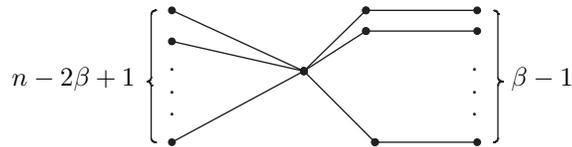}
  \caption{ $T_{n, \beta}$ }
\end{center}
\end{figure}

Let $T_{n, \beta}$ be the tree obtained from the star graph $S_{n-\beta+1}$ by attaching a pendant edge to each of certain $\beta-1$ non-central vertices of $S_{n-\beta+1}$; see Fig. 4. It is obvious to see that $T_{n, \beta}$ contains an $\beta$-matching. In the following, we show that $T_{n, \beta}$ is the tree with the maximum $\xi^{ee}$ among $n$-vertex trees with matching number $\beta$. Let $P_t(a,b)$ be the $n$-vertex graph obtained by attaching $a$ and $b$ leaves to the endvertices of $P_t,$ respectively. Clearly, $a+b+t=n.$

\begin{thm}
Let $T$ be an $n$-vertex ($n\geqslant2\beta$) tree with matching number $\beta$.
 \begin{wst}
\item[{\rm (i)}]If $\beta=1,$ then $\xi^{ee}(T)=\xi^{ee}(S_n)=\frac{3n-3}{2}$.
\item[{\rm (ii)}]If $\beta=2$,  then $\xi^{ee}(T)\leqslant \frac{5n-4}{6}$ with equality if and only if $T\cong P_2(a, b),$ where $a+b+2=n$.
\item[{\rm (iii)}]If $\beta \geqslant 3$, then $\xi^{ee}(T)\leqslant \frac{10n-3\beta-7}{12}$ with equality if and only if $T\cong T_{n, \beta}$, where $T_{n, \beta}$ is depicted in Fig. 4.
\end{wst}
\end{thm}
\begin{proof}\ (i)\ If $\beta=1$, then there is just one such $n$-vertex tree, $S_n$. By a direct calculation, we have $\xi^{ee}(T)=\frac{3n-3}{2}$.

(ii)\ If $\beta=2$, then $T\cong P_2(a,b)$ with $a+b+2=n$, or $T\cong P_3(s,t)$ with $s+t+3=n$. By a simple computing, we have $\xi^{ee}(T)\leqslant \frac{5n-4}{6}$, and the equality holds if and only if $T\cong P_2(a,b)$ with $a+b+2=n$.

(iii)\ Choose an $n$-vertex tree $T$ with matching number $\beta \geqslant 3$ such that $\xi^{ee}(T)$ is as large as possible. If there is a pendant path $v_0v_1v_2\ldots v_{l-3}v_{l-2}v_{l-1}v_l$ attached at vertex $v_0$ in $T$ with $l>2$, then we let $T_1=T-v_{l-2}v_{l-1}+v_0v_{l-1}.$ It is routine to check that $\beta(T_1)=\beta(T)$. By Corollary 3.2, we have $\xi^{ee}(T)<\xi^{ee}(T_1)$, a contradiction.  Therefore, the length of each of the pendant paths in $T$ is one or two.
\begin{figure}[h!]
\begin{center}
\psfrag{1}{$w$}
\psfrag{2}{$G$}
\psfrag{3}{$p$}
\psfrag{4}{$q+1$}
\psfrag{5}{$v$}
\psfrag{6}{$T'$}
\psfrag{7}{$q$}
\psfrag{8}{$T$}
\psfrag{a}{$p+1$}
\psfrag{b}{$q-1$}
\psfrag{c}{$T''$}\psfrag{d}{$\Rightarrow$}\psfrag{r}{$\Leftarrow$}
\includegraphics[width=150mm]{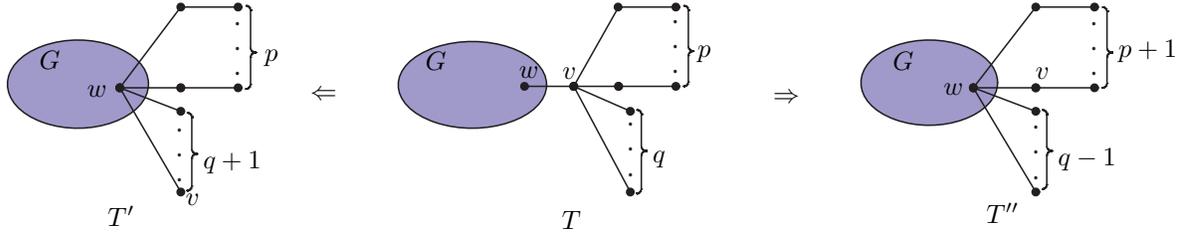}
  \caption{ $T\Rightarrow T'$ and $T\Rightarrow T''$ }
\end{center}
\end{figure}
Hence, we may assume that there are $p$ $P_3$'s and $q$  $P_2$'s attached at $v$; see Fig. 5.

First we consider that $w$ is a parent of $v$ or both $w$ and $v$ are two adjacent center vertices of $T$. We proceed by considering the following two possible cases. 

If the vertex $w$ is not perfectly matched in $T$, then there exists a matching $M$ of maximum cardinality such that $w$ is not $M$-matched. Applying $\alpha$-transformation to $T$ yields the tree $T'$ (see $T\Rightarrow T'$ in Fig. 5). Note that $\beta(T')=\beta(G-w)+p+1=\beta(G)+p+1=\beta(T)$. That is to say, $T'$ is an $n$-vertex tree with matching number $\beta.$ By Theorem 3.3, $\xi^{ee}(T)<\xi^{ee}(T')$, a contradiction.

If the vertex $w$ is perfectly matched in $T$, then for each matching $M$ with maximum cardinality, $w$ is $M$-matched. Applying $\rho$-transformation to $T$ yields the tree $T''$ (see $T\Rightarrow T''$ in Fig. 5).  Obviously, $\beta(T'')=\beta(G)+p+1=\beta(T)$. Based on Theorem 3.1, we have $\xi^{ee}(T)<\xi^{ee}(T''),$ a contradiction.

Now we consider that $v$ is a parent of $w.$ Then the subgraph $G$ is either a star or is an isolated vertex $w.$ If the former happens, then by a similar discussion as above, we may obtain another $n$-vertex tree $T'$ with matching number $\beta$ satisfying $\xi^{ee}(T)<\xi^{ee}(T')$, a contradiction. Hence, we obtain that $G\cong K_1.$  Hence, we obtain that, in $T$, $p=\beta-1$ and $q=n-2\beta+1$, i.e., $T\cong T_{n, \beta}$, which is the unique tree with maximum $\xi^{ee}$ among $n$-vertex trees with matching number $\beta\geqslant3$. By a simple calculation, we have that $\xi^{ee}(T_{n, \beta})=\frac{10n-3\beta-7}{12}.$

This completes the proof.
\end{proof}
Recall that $\mathscr{D}_n^\gamma$ denote the set of all $n$-vertex trees with domination number $\gamma$.
\begin{lem}[Haynes et. al. 1998]
For a graph $G$, we have $\gamma(G)\leqslant \beta(G)$.
\end{lem}
\begin{lem}
If $T\in \mathscr{D}_n^\gamma$ has the maximum $\xi^{ee}$, then $\gamma(T)=\beta(T)=\gamma$.
\end{lem}
\begin{figure}[h!]\label{fig.6}
\begin{center}
  \psfrag{1}{$v_1$}\psfrag{2}{$v_{1'}$}\psfrag{3}{$v_2$}\psfrag{4}{$v_{2'}$}\psfrag{5}{$v_\gamma$}
  \psfrag{6}{$v_{\gamma'}$}\psfrag{7}{$w_1$}\psfrag{8}{$w_2$}\psfrag{a}{$T$}\psfrag{b}{$T'$}
  \includegraphics[width=100mm]{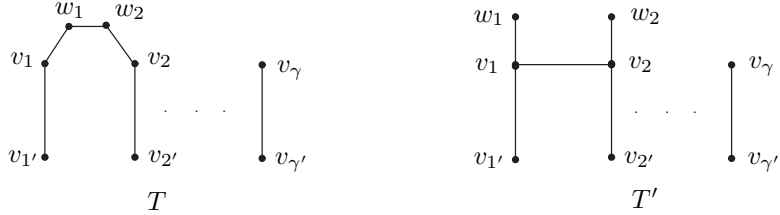}\\
  \caption{The structures of $T$ and $T'.$ }
\end{center}
\end{figure}
\begin{proof}
It suffices to show that $\beta(T)\leqslant\gamma(T)$ by Lemma 4.3. Otherwise, we have $\gamma=\gamma(T)<\beta(T)$. Assume that $S=\{v_1, v_2, \ldots, v_{\gamma}\}$ is a minimum dominating set of $T$. Then there exist $\gamma$ independent edges, say $v_1v_1', v_2v_2', \ldots, v_{\gamma}v_{\gamma}'$. Let $M'=\{v_iv_i': i=1, 2, \ldots, \gamma\}$, which is obvious a matching of $T$.

If $M'$ is contained in a maximum matching, say $M$, of $T$, then there must exist another edge, say $w_1w_2$, with no shared endvertex with each of $v_iv_i',\, i=1,2,\ldots, \gamma$ (based on $\gamma(T)<\beta(T)$).

If $M'$ is not contained in any maximum matching of $T$, then $M'$ is a maximal matching. Note that $\gamma(T)<\beta(T)$; hence there exists an $M'$-augmenting path of length $2t+1,\, t\leqslant \gamma$ and denote it by $P=u_1u_2\ldots u_{2t+1}u_{2t+2}$, where $u_{2k}u_{2k+1}=v_kv_k'$, $k=1, 2, \ldots, t$. Thus, for each adjacent pair $u_i, u_j\in S\cap V_P$ with $i<j$, we have $2\leqslant j-i\leqslant3$. In particular,  there must exist a pair of such vertices $u_i, u_j\in S\cap V_P$ satisfying $j-i=3$ since the first (resp. last) domination vertex is $u_2$ (resp. $u_{2t+1}$). Suppose $i_0$ is the smallest value of $i$ for which $u_{i_0}, u_{i_0+3}\in S\cap V_P$. Let
$$
  M''=\{u_1u_2, u_3u_4, \ldots, u_{i_0-1}u_{i_0}\}\cup\{u_{i_0+3}u_{i_0+4}, u_{i_0+5}u_{i_0+6}, \ldots, u_{2t+1}u_{2t+2}\}\cup \{v_{t+1}v'_{t+1},\ldots, v_{\gamma}v_{\gamma}'\}.
$$
Then $M''$ is a matching of cardinality $\gamma$, and each vertex in $S$ is $M''$-matched. It is routine to check that the edge $u_{i_0+1}u_{i_0+2}$ is independent of each edge from $M''$. Hence, in this case we also obtain $\gamma+1$ independent edges $v_1v_1', v_2v_2', \ldots, v_{\gamma}v_{\gamma}', u_{i_0+1}u_{i_0+2}$ in $T$.

Summarize the discussion as above, we conclude that there exists $\gamma+1$ independent edges $v_1v_1', v_2v_2', \ldots, v_{\gamma}v_{\gamma}',\linebreak w_1w_2$ in $T$.

In what follows we show that $w_1, w_2$ are dominated by two different vertices from $S$. In fact, if this is not true, then there will occur a triangle in $T$, which is a contradiction to the fact that $T$ is a tree. Without loss of generality, assume that $w_1, w_2$ are dominated by $v_1, v_2$, respectively.  By $\alpha$-transformation of $T$ on the edges $v_1w_1$ and $v_2w_2$, we can obtain a new tree $T'\in \mathscr{D}_n^\gamma$ such that $\xi^{ee}(T')>\xi^{ee}(T)$ by Theorem 3.3 (see Fig. 6), a contradiction. That is, $\beta(G)\leqslant\gamma(G)$. Together with Lemma 4.3, we obtain $\gamma(T)=\beta(T)=\gamma$, as desired.
\end{proof}

From Theorem 4.2 and Lemma 4.4, we can easily get the next result.
\begin{thm}
For any tree $T\in \mathscr{D}_n^\gamma$.
 \begin{wst}
\item[{\rm (i)}]If $\gamma=1$, then $\xi^{ee}(T)=\xi^{ee}(S_n)=\frac{3n-3}{2}$.
\item[{\rm (ii)}]If $\gamma=2$, then $\xi^{ee}(T)\leqslant \frac{5n-4}{6}$ with equality if and only if $T\cong P_2(a, b),$ where $a+b+2=n$.
\item[{\rm (iii)}]If $\gamma\geqslant3$, then $\xi^{ee}(T)\leqslant \frac{10n-3\gamma-7}{12}$ with equality if and only if $T\cong T_{n, \gamma}$.
\end{wst}
\end{thm}

Let $\mathscr{J}_n^{p,q}$ be the set of all $n$-vertex trees with a ($p, q$)-bipartition, $q\geqslant p\geqslant 1$. Note that $\mathscr{J}_n^{1, n-1}$ contains just $S_n$; $\mathscr{J}_n^{2, n-2}=\{P_3(a, b)| a+b=n-3\}$, where $P_3(a,b)$ is an $n$-vertex trees obtained by attaching $a$ and $b$ leaves to the endvertices of $P_3$, respectively. By a direct calculation, we have $\xi^{ee}(P_3(a, b))=\frac{7n-1}{12}$ for $a\geqslant1$ and $b\geqslant1$, which is obvious if $a\geqslant2$ and $b\geqslant2$, $\xi^{ee}(P_3(a, b))=\xi^{ee}(P_3(a-1, b+1))=\xi^{ee}(P_3(a+1, b-1))$. Furthermore,
\begin{thm}
Given positive integers $p, q$ with $q\geqslant p>2$ and $p+q=n$, then, for $T\in \mathscr{J}_n^{p, q}$, one has $\xi^{ee}(T)\leqslant \frac{5n-4}{6}$ with  equality if and only if $T\cong P_2(p,q)$.
\end{thm}
\begin{proof}
For a $T\in \mathscr{J}_n^{p, q}$ with $q\geqslant p>2$ and $p+q=n$, repeatedly applying $\theta$-transformation to $T$ yields that $P_2(p,q)$ is the unique tree in $\mathscr{J}_n^{p, q}$ such that it has the maximum $\xi^{ee}$. By an elementary calculation, we have $\xi^{ee}(T)\leqslant \frac{5n-4}{6}$.
\end{proof}

Let $\mathscr{H}_n^d$ be the set of $n$-vertex trees of diameter $d$. For any graph $T$ in $\mathscr{H}_n^d$, let $P_{d+1}=v_0v_1\ldots v_d$ be a diametric path of length $d$ in $T$. In particular, let $C(a_1, a_2, \ldots, a_{d-1})$ be an $n$-vertex tree obtained from $P_{d+1}$ by attaching $a_i$ pendant edges to vertex  $v_i$, $i=1, 2, \dots, d-1$. Obviously, $\diam(C(a_1, \ldots, a_{d-1}))=d$ and $n=d+1+\sum_{i=1}^{d-1}a_i$. Denote $C_{n ,d}=C(0, \ldots, 0, a_{\frac{d}{2}}, 0, \ldots, 0)$ if $d$ is even and  $C_{n ,d}=C(0, \ldots, 0, a_{\frac{d-1}{2}}, a_{\frac{d+1}{2}}, 0, \ldots, 0)$ if $d$ is odd.

Yu et al. \cite{VY} showed that $C_{n, d}$ is the unique graph having the maximal $\xi^{ee}$ value among $\mathscr{H}_n^d$. Here, we determine sharp upper bounds on $\xi^{ee}(T)$ of graphs among $\mathscr{H}_n^d\setminus \{C_{n, d}\}$.
\begin{figure}[h!]
\begin{center}
\psfrag{1}{$v_0$}
\psfrag{2}{$v_1$}
\psfrag{3}{$v_{\frac{d}{2}-1}$}
\psfrag{4}{$v_{\frac{d}{2}}$}
\psfrag{5}{$v_{\frac{d}{2}+1}$}
\psfrag{6}{$v_{d-1}$}
\psfrag{7}{$v_d$}
\psfrag{8}{$v_{\frac{d-1}{2}-1}$}
\psfrag{9}{$v_{\frac{d-1}{2}}$}
\psfrag{e}{$v_{\frac{d+1}{2}}$}
\psfrag{f}{$v_{\frac{d+1}{2}+1}$}
\psfrag{c}{$p$}\psfrag{H}{$T^1$}\psfrag{I}{$T^2$}
\psfrag{d}{$q$}\psfrag{J}{$T^3$}\psfrag{K}{$T^4$}
\psfrag{s}{$s$}\psfrag{L}{$T^5$}\psfrag{M}{$T^6$}
\psfrag{t}{$t$}\psfrag{N}{$T^7$}
\psfrag{a}{$n-d-2$}
\psfrag{x}{$p+q=n-d-2$}
\psfrag{b}{$n-d-3$}
\psfrag{k}{$\mathcal{F}_{n, d}$}
\psfrag{o}{$\mathcal{F}_{n, d}^{s, t}$}
\psfrag{g}{$s+t=n-d-3$}
\includegraphics[width=150mm]{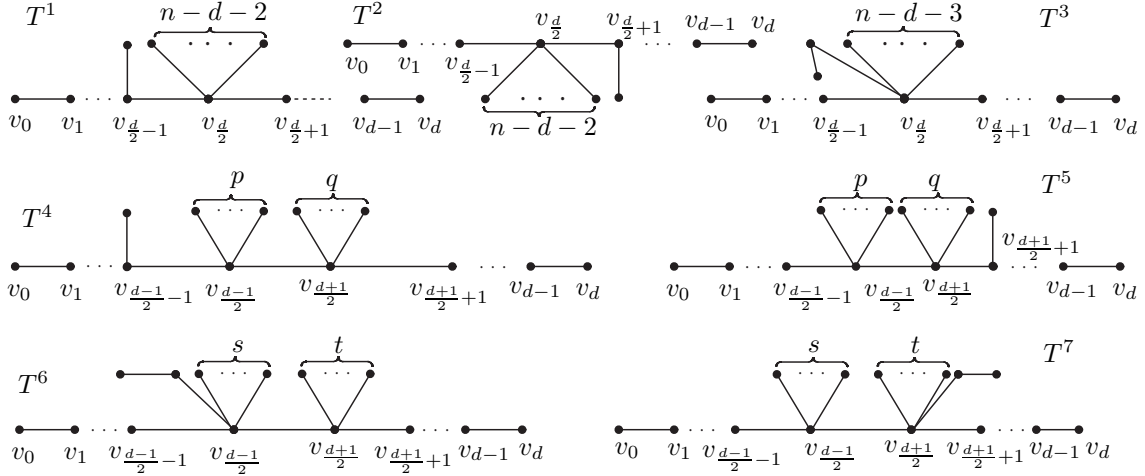}
  \caption{Graphs $T^1, T^2, T^3, T^4, T^5, T^6$ and $T^7$, where $p+q=n-d-2$ and $s+t=n-d-3$.}
\end{center}
\end{figure}
\begin{thm}
Let $T\in \mathscr{H}_n^d\setminus \{C_{n, d}\}$  with $d\geqslant4$. Then
 $$  \xi^{ee}(T)\leqslant\left\{
             \begin{array}{ll}
              \sum_{i=0}^{\frac{d}{2}-2}\frac{4}{d-i}+\frac{2n-2d-2}{d}+\frac{2n-2d+6}{d+2}+\frac{2}{d+4}, & \hbox{if $d$ is even,}\\[5pt]
               \sum_{i=0}^{\frac{d-1}{2}-2}\frac{4}{d-i}-\frac{2}{d}+\frac{2n-2d+4}{d+1}+\frac{2n-2d+6}{d+3}+\frac{2}{d+5}, & \hbox{if $d$ is odd.}
             \end{array}
           \right.
$$
The equality holds if and only if $T\cong T^1, T^2$ or $T^3$ for even $d$ and $T\cong T^4, T^5, T^6$ or $T^7$ for odd $d$, where $T^1, T^2, \ldots, T^7$ are depicted in Fig. 7.
\end{thm}
\begin{proof}
Choose an $n$-vertex tree $T$ in $\mathscr{H}_n^d\setminus \{C_{n, d}\}$ such that $\xi^{ee}(T)$ is as large as possible. Denote the diametric path of $T$ by $P=v_0v_1\ldots v_d$.

First, we consider even $d$. In this case, let $N_2(v_{\frac{d}{2}})=\{x\in V_T\setminus V_P: d_T(x, v_{\frac{d}{2}})=2\}$.

If $N_2(v_{\frac{d}{2}})=\emptyset$, then all vertices of $V_T\setminus V_P$ must be pendant vertices attached at some non-pendant vertices of $P$, otherwise by Theorem 3.3 we may get another $n$-vertex tree $T_1$ in $\mathscr{H}_n^d\setminus \{C_{n, d}\}$ such that $\xi^{ee}(T)<\xi^{ee}(T_1)$, a contradiction. Assume that the non-pendant vertex $v_i$ ($i\neq \frac{d}{2}$) on $P$ is adjacent to some pendant vertices, say $u_1, u_2, \ldots,  u_t$, $t\geqslant2$. Let
$$
   T'=T-\{v_iu_2, v_iu_3, \ldots, v_iu_t\}+\{v_{\frac{d}{2}}u_2, v_{\frac{d}{2}}u_3, \ldots, v_{\frac{d}{2}}u_t\}.
$$
Thus, we have
\begin{eqnarray*}
\xi^{ee}(T)-\xi^{ee}(T')&=&\sum_{k=2}^{t}\left(\frac{1}{\varepsilon_T(v_i)}+\frac{1}{\varepsilon_T(u_k)}\right)
-\sum_{k=2}^{t}\left(\frac{1}{\varepsilon_{T'}(v_{\frac{d}{2}})}+\frac{1}{\varepsilon_{T'}(u_k)}\right)\\
&=&\frac{t-1}{\varepsilon_T(v_i)}+\frac{t-1}{\varepsilon_T(u_k)}-\frac{t-1}{\varepsilon_{T'}(v_{\frac{d}{2}})}-\frac{t-1}{\varepsilon_{T'}(u_k)}\\
&=&\frac{t-1}{\varepsilon_T(v_i)}+\frac{t-1}{\varepsilon_T(v_i)+1}-\frac{t-1}{\varepsilon_{T'}(v_{\frac{d}{2}})}-\frac{t-1}{\varepsilon_{T'}(v_{\frac{d}{2}})+1}.
\end{eqnarray*}
Note that $\varepsilon_T(v_i)>\frac{d}{2}=\varepsilon_{T'}(v_{\frac{d}{2}});$ hence $\xi^{ee}(T)-\xi^{ee}(T')<0,$ a contradiction. Therefore, every non-pendant vertex $v_i$ ($i\neq \frac{d}{2}$) on $P$ is adjacent to at most one pendant vertex. It is routine to check that there is at least one vertex in $\{v_1,v_2,\ldots, v_{\frac{d}{2}-1}, v_{\frac{d}{2}+1}, \ldots, v_{d-1}\}$ being adjacent to just one pendant vertex; otherwise $T\cong C_{n,d},$ a contradiction. So in what follows, we show that there exists a $x$ in $PV(T)$ such that it is only adjacent to $v_{\frac{d}{2}-1}$ or $v_{\frac{d}{2}+1}.$
If this is not true, then either the pendant vertex $x\sim v_i$ with $i\not= 0,\frac{d}{2}-1,\frac{d}{2}+1$ and $d,$ or one pendant vertex $x\sim v_{\frac{d}{2}-1}$ and another pendant vertex $y\sim v_{\frac{d}{2}+1}$. If the former happens, then let
$$  T'=\left\{
             \begin{array}{ll}
              T-v_ix+v_{i+1}x, & \hbox{if $1\leqslant i <\frac{d}{2}-1,$ }\\[5pt]
              T-v_ix+v_{i-1}x, & \hbox{if $\frac{d}{2}+1 <i\leqslant d-1.$}
             \end{array}
           \right.
$$
If the latter happens, we let $T^1=T-v_{\frac{d}{2}-1}x+v_{\frac{d}{2}}x$ or let $T^2=T-v_{\frac{d}{2}+1}y+v_{\frac{d}{2}}y$. By Theorem 3.1, we obtain
$\xi^{ee}(T)<\xi^{ee}(T'),\, \xi^{ee}(T)<\xi^{ee}(T^1)$ and $\xi^{ee}(T)<\xi^{ee}(T^2)$, a contradiction. Hence, we obtain that $T\cong T^1$ or $T^2$ if $N_2(v_{\frac{d}{2}})=\emptyset$, where $T^1, T^2$ are depicted in Fig. 7.

If $N_2(v_{\frac{d}{2}})\neq\emptyset$, then by Theorems 3.3 and 3.6, we obtain that $V_T\setminus V_P\subseteq PV(T).$ Furthermore, we may partition $V_T\setminus V_P$ as $V_1\cup V_2$, where $V_1= N_2(v_{\frac{d}{2}})$ and each vertex in $V_2$ is adjacent to some vertex on the diametric path $P.$

If $N_2(v_{\frac{d}{2}})=\{w_1, w_2, \ldots, w_k\}$ with $k\geqslant 2$, then let $T_1=T-\{w_2, w_3, \ldots, w_k\}+\{v_{\frac{d}{2}}w_2, v_{\frac{d}{2}}w_3, v_{\frac{d}{2}}w_k\}.$ By a similar discussion as in the proof of Theorem 3.6, we may obtain that $\xi^{ee}(T)<\xi^{ee}(T_1)$, a contradiction. Hence, we have $|N_2(v_{\frac{d}{2}})|= 1.$

If $V_2$ contains a vertex, say $x$, such that $x$ is adjacent to $v_i$ on $P$ with $i\not= 0,\frac{d}{2}$ and $d,$ then let
$$  T'=\left\{
             \begin{array}{ll}
              T-v_ix+v_{i+1}x, & \hbox{if $1\leqslant i \leqslant\frac{d}{2}-1,$ }\\[5pt]
              T-v_ix+v_{i-1}x, & \hbox{if $\frac{d}{2}+1 \leqslant i\leqslant d-1.$}
             \end{array}
           \right.
$$
By Theorem 3.1, we obtain $\xi^{ee}(T)<\xi^{ee}(T')$, a contradiction. That is $T\cong T^3$ (see Fig. 7). By direct computing, we have
$$
\xi^{ee}(T^1)=\xi^{ee}(T^2)=\xi^{ee}(T^3)=\sum_{i=0}^{\frac{d}{2}-2}\frac{4}{d-i}+\frac{2n-2d-2}{d}+\frac{2n-2d+6}{d+2}+\frac{2}{d+4},
$$
as desired.

If $d$ is odd, by a similar discussion as in the proof for even $d$, we may also show our result is true. We omit the procedure here.
\end{proof}

Let $\mathscr{T}_{n,d}^{p, q}$ be the set of $n$-vertex trees obtained from $P_{d+1}=v_0v_1\ldots v_d$ by attaching $p$ and $q$ pendant vertices at $v_1$ and $v_{d-1}$, respectively, where $p+q=n-d-1$. Yu et al. \cite{VY} showed that each tree in $\mathscr{T}_{n,d}^{p, q}$ having the minimal $\xi^{ee}$-value among all the set of $n$-vertex trees of diameter $d$. In the rest of this section, we characterize all the trees $\mathscr{H}_n^d\setminus \mathscr{T}_{n,d}^{p, q}$ with the minimal $\xi^{ee}$-value.
\begin{figure}[h!]
\begin{center}
\psfrag{1}{$v_0$}\psfrag{A}{$T'$}
\psfrag{2}{$v_1$}\psfrag{B}{$T''$}
\psfrag{3}{$v_2$}\psfrag{C}{$\hat{T}$}
\psfrag{4}{$v_3$}\psfrag{r}{$r$}
\psfrag{5}{$v_{d-2}$}\psfrag{s}{$s$}
\psfrag{6}{$v_{d-1}$}\psfrag{t}{$t$}
\psfrag{7}{$v_d$}\psfrag{u}{$h$}
\psfrag{8}{$p$}
\psfrag{9}{$q$}
\includegraphics[width=150mm]{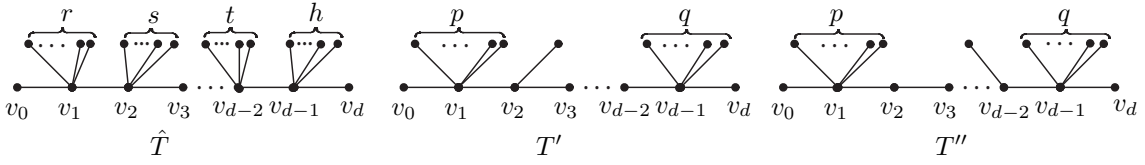}
  \caption{Graph $\hat{T}$ with $r+s+t+h=n-d-1, s+t>0$ and graphs $T', T''$ with $p+q=n-d-2.$}
\end{center}
\end{figure}
\begin{thm}
Let $T\in \mathscr{H}_n^d\setminus \mathscr{T}_{n,d}^{p, q}$. Then
$$  \xi^{ee}(T)\geqslant\left\{
             \begin{array}{ll}
              \sum_{i=1}^{\frac{d}{2}-1}\frac{4}{d-i}+\frac{n-d+4}{d}+\frac{n-d-1}{d-1}+\frac{1}{d-2}, & \hbox{if $d$ is even,}\\[5pt]
              \sum_{i=1}^{\frac{d-1}{2}}\frac{4}{d-i}+\frac{n-d}{d}+\frac{n-d-1}{d-1}+\frac{1}{d-2}, & \hbox{if $d$ is odd.}
             \end{array}
           \right.
$$
The equality holds if and only if $T\cong T'$ or $T''$, where $T'$ and $T''$ are depicted in Fig. 8.
\end{thm}
\begin{proof}
Let $T\in \mathscr{H}_n^d\setminus \mathscr{T}_{n,d}^{p, q}$ be a tree with diameter $P=v_0v_1\ldots v_d$. Denote the component in $T-\{v_iv_{i-1}, v_iv_{i+1}\}$ containing $v_i$ by $T_i$, $i=1,2,\ldots, d-1.$ That is, $T$ is obtained from $P=v_0v_1\ldots v_d$ by attaching $T_i$ to $v_i$ for $i=1,2,\ldots, d-1.$ Let $V^0=\{x\in V_T\setminus (V_P\cup V_{T_1}\cup V_{T_{d-1}}): \varepsilon_T(x)=d\}.$ It is routine to check that $V^0\subseteq PV(T)$. Put
$$
T^*=T-V^0+\{v_1x: x\in V^0\}.
$$
Clearly, $\xi^{ee}(T^*)=\xi^{ee}(T).$ It is easy to see that $T^*\in \mathscr{H}_n^d\setminus \mathscr{T}_{n,d}^{p, q}$ and there is at least one pendant vertex, say $w,$ such that $w\not\sim v_1$ and $w\not\sim v_{d-1}.$ Note that $P=v_0v_1\ldots v_d$ is as well a diametric path of $T^*.$ Let $V^1=N_{T^*}(v_1)\cup N_{T^*}(v_2)\cup N_{T^*}(v_{d-2})\cup N_{T^*}(v_{d-1})\cup V_P.$ If $V_{T^*}\setminus V^1=\emptyset,$ then $T^*$ has the same structure as that of $\hat{T}$ as depicted in Fig. 8; If $V_{T^*}\setminus V^1\not=\emptyset,$ then put
$$
  T^{**}=T^*-V_{T^*}\setminus V^1+\{v_2x: x\in V_{T^*}\setminus V^1\}.
$$
It is routine to check that $T^{**}\in \mathscr{H}_n^d\setminus \mathscr{T}_{n,d}^{p, q}$ and $\varepsilon_{T^*}(x)\leqslant \varepsilon_{T^{**}}(x)=d-1,\, d_{T^{**}}(x) \leqslant d_{T^*}(x)$ for all $x\in V_{T^*}\setminus V^1.$ In particular, there exists at least one vertex, say $a\in V_{T^*}\setminus V^1$, satisfying $\varepsilon_{T^*}(a)< \varepsilon_{T^{**}}(a)=d-1,\, d_{T^{**}}(a) \leqslant d_{T^*}(a)$. Thus, by a direct calculation we may get $\xi^{ee}(T^{**})<\xi^{ee}(T^*)$ and $T^{**}$ has the same structure as that of $\hat{T}$ as depicted in Fig. 8

Hence, in what follows we are to determine the extremal graph from $\hat{T}$ as depicted in Fig. 8. Note that in $\hat{T}$, $s+t>0.$ Without loss of generality, we assume that $s>0$. If $t>0$, then we move all the $t$ pendant edges from $v_{d-2}$ to $v_2$ and denote the resultant graph by $\hat{T}'$. It is easy to see that $\xi^{ee}(\hat{T}')=\xi^{ee}(\hat{T})$. So we assume that $t=0,s\geqslant 1$ in $\hat{T}$. If $s=1$, then we have $\hat{T}\cong T'$ or $T''$ (see Fig. 8). If $s\geqslant 2$, then let $\bar{T}$ be the tree obtained from $\hat{T}$ by moving $s-1$ pendants from $v_2$ to $v_1.$ By a direct computing, we have
\begin{eqnarray*}
  \xi^{ee}(\bar{T})-\xi^{ee}(\hat{T}) &=& \frac{d_{\hat{T}}(v_1)+s-1}{\varepsilon_{\hat{T}}(v_1)}+\frac{3}{\varepsilon_{\hat{T}}(v_2)}+\frac{s-1}{\varepsilon_{\hat{T}}(v_1)+1}-\left( \frac{d_{\hat{T}}(v_1)}{\varepsilon_{\hat{T}}(v_1)}+\frac{s+2}{\varepsilon_{\hat{T}}(v_2)}+\frac{s-1}{\varepsilon_{\hat{T}}(v_2)+1}\right) \\
   &=& (s-1)\left(\frac{1}{\varepsilon_{\hat{T}}(v_1)}-\frac{1}{\varepsilon_{\hat{T}}(v_2)}+\frac{1}{\varepsilon_{\hat{T}}(v_1)+1}-\frac{1}{\varepsilon_{\hat{T}}(v_2)+1} \right) \\
   &<&0.
\end{eqnarray*}
The last inequality follows by $s\geqslant 2$ and $\varepsilon_{\hat{T}}(v_1)=d-1>\varepsilon_{\hat{T}}(v_2)=d-2.$ Hence, we obtained that the graph $\hat{T}\cong T'$ or $T''$ having the minimal $\xi^{ee}$-value among $\mathscr{H}_n^d\setminus \mathscr{T}_{n,d}^{p, q}$, where $T'$ and $T''$ are depicted in Fig. 8. By direct calculation, we have
$$  \xi^{ee}(\hat{T})=\xi^{ee}(T')=\xi^{ee}(T'')\geqslant\left\{
             \begin{array}{ll}
              \sum_{i=1}^{\frac{d}{2}-1}\frac{4}{d-i}+\frac{n-d+4}{d}+\frac{n-d-1}{d-1}+\frac{1}{d-2}, & \hbox{if $d$ is even,}\\[5pt]
              \sum_{i=1}^{\frac{d-1}{2}}\frac{4}{d-i}+\frac{n-d}{d}+\frac{n-d-1}{d-1}+\frac{1}{d-2}, & \hbox{if $d$ is odd,}
             \end{array}
           \right.
$$
as desired.
\end{proof}
\section{\normalsize Summary and conclusions}
In this article, we studied the the total reciprocal edge-eccentricity of graphs, which was introduced by Gupta, Singh and Madan in \cite{4} and derive some monotonicity properties on this novel graph invariant under some edge-graph transformations. In view of \cite{11,VY}, there is a lack of further analytical results in the scientific literature when studying this distance-degree-based graph invariants on trees. We obtained some sharp bounds on the total reciprocal edge-eccentricity of trees in terms of graph parameters such as pendants, matching number, domination number, diameter, vertex bipartition, et al,  which extended some of the results obtained in \cite{11,VY}. As a future work, we want to explore general methods
to show the extremal values of the REE for characterizing the structural properties of graphs.

\end{document}